\documentclass[11pt]{article} 
\usepackage{amssymb,amsfonts,amsmath,amsthm,amscd,stackrel,dsfont,mathrsfs,eqsection}
\RequirePackage[colorlinks,citecolor=blue,urlcolor=blue]{hyperref}

\newcommand{\eqnsection}{\renewcommand{\theequation}{\thesection.\arabic{equation}}
      \makeatletter \csname @addtoreset\endcsname{equation}{section}\makeatother}

\newtheorem{theorem}{Theorem}[section]

\newtheorem{lemma}[theorem]{Lemma}
\newtheorem{proposition}[theorem]{Proposition}

\newtheorem{remark}[theorem]{Remark}

\newcommand{\be}{\begin{equation}}
\newcommand{\ee}{\end{equation}}
\newcommand{\beq}{\begin{equation*}}
\newcommand{\eeq}{\end{equation*}}
\newcommand{\bq}{\begin{eqnarray}}
\newcommand{\eq}{\end{eqnarray}}

\newcommand\ceil[1]{\lceil#1\rceil}

\def\P{{\mathbb{P}}}

\newcommand{\eps}{\varepsilon}

\newcommand{\EE}{\mathbb E}

\newcommand{\RR}{\mathbb R}
\def\reals{\RR}

\newcommand{\PP}{\mathbb P}

\newcommand{\cK}{{\mathcal K}}

\newcommand{\cN}{{\mathcal N}}

\newcommand{\abbr}[1]{{\sc\lowercase{#1}}}

\begin{document}

\title{Persistence of iterated partial sums}
\author{{\sc Amir Dembo}\thanks{The research was supported in part by NSF grant DMS-0806211.}
\and {\sc Fuchang Gao}\thanks{This research was partially done when the second named author visited Stanford university in March 2010. The authors also thank NSF grant DMS-0806211 for the support.}}
\date{\today}
\maketitle
\begin{abstract}
Let $S_n^{(2)}$ denote the iterated partial sums. That is,
$S_n^{(2)}=S_1+S_2+\cdots +S_n$, where $S_i=X_1+X_2+\cdots+X_i$.
Assuming $X_1, X_2,\ldots,X_n$ are integrable, zero-mean,
i.i.d. random variables, we show that the persistence probabilities
$$p_n^{(2)}:=\PP\left(\max_{1\le i \le n}S_i^{(2)}<  0\right) \le c\sqrt{\frac{\EE|S_{n+1}|}{(n+1)\EE|X_1|}},$$
with $c \le 6 \sqrt{30}$ (and $c=2$ whenever $X_1$ is symmetric).
Furthermore, the converse inequality holds whenever
$\PP(-X_1>t)\asymp e^{-\alpha t}$ for some $\alpha>0$
or $\PP(-X_1>t)^{1/t} \to 0$ as $t \to \infty$.
Consequently, for these random variables we
have that $p_n^{(2)}\asymp n^{-1/4}$ if $X_1$ has finite second moment.
In contrast, we show that for any $0 < \gamma < 1/4$ there exist
integrable, zero-mean random variables for which the rate of
decay of $p_n^{(2)}$ is $n^{-\gamma}$.
\end{abstract}

\section{Introduction}
The estimation of probabilities of rare events is one of the
central themes of research in the theory of probability. Of
particular note are \emph{persistence} probabilities, formulated as
\bq
q_n = \PP\Big(\max_{1\le k\le n} Y_k < y \Big),\label{max}
\eq
where $\{Y_k\}_{k=1}^n$ is a sequence of zero-mean random variables.
For independent $Y_i$ the persistence probability is easily
determined to be the product of $\PP(Y_k < y)$ and to a large
extent this extends to the case of sufficiently weakly dependent
and similarly distributed $Y_i$, where typically $q_n$ decays
exponentially in $n$. In contrast, in the classical case of partial sums,
namely $Y_k=S_k=\sum_{i=1}^k X_i$ with $\{X_j\}$ i.i.d.
zero-mean random variables, it is well known that $q_n = O(n^{-1/2})$
decays as a power law. This seems to be one of the very few cases
in which a power law decay for $q_n$ can be proved and its exponent
is explicitly known. Indeed, within the large class of similar
problems where dependence between $Y_i$ is strong enough to rule out
exponential decay, the behavior of $q_n$ is very sensitive to
the precise structure of dependence between the variables
$Y_i$ and even merely determining its asymptotic rate can be very
challenging (for example, see \cite{DPSZ} for recent results in case
$Y_k = \sum_{i=1}^n X_i (1-c_{k,n})^i$ are the values of a random
Kac polynomials evaluated at certain non-random $\{c_{k,n}\}$).

We focus here on iterated sums of i.i.d. squared integrable, zero-mean,
random variables $\{X_i\}$. That is, with $S_n=\sum_{k=1}^n X_k$ and
\bq
S^{(2)}_n=\sum_{k=1}^n S_k=\sum_{i=1}^n(n-i+1)X_i\,,
\label{S2}
\eq
we are interested in the asymptotics as $n \to \infty$
of the persistence probabilities
\bq
p_n^{(2)}(y):=\P\left(\max_{1\le k\le n}S_k^{(2)}<
y\right),\  \overline{p}_n^{(2)}(y):=\P\left(\max_{1\le k\le n}S_k^{(2)}\le
y\right)\,, \label{prob}
\eq
where $y \ge 0$ is independent of $n$.
With $y \ll n$ it immediately
follows from Lindeberg's \abbr{clt} that $p_n^{(2)}(y)\to 0$
as $n\to\infty$ and our goal is thus to find a sharp rate for
this decay to zero.

Note that for any fixed $y>0$ we have that
$\overline{p}_n^{(2)} (y) \sim p_n^{(2)}(y) \sim p_n^{(2)}(0)$
up to a constant depending only on $y$. Indeed, because
$\EE X^->0$, clearly $\PP(X_1<-\eps)>0$ for $\eps=y/k$
and some integer $k \ge 1$. Now, for any $n\ge 1$ and $z \ge 0$,
$$
\overline{p}^{(2)}_n(z) \ge
p_n^{(2)}(z) \ge \PP(X_n<-\eps)\overline{p}_{n-1}^{(2)}(z+\eps)\ge
\PP(X_1<-\eps) \overline{p}_{n}^{(2)}(z+\eps)
$$
and applying this inequality for $z=i \eps$, $i=0,1,\ldots,k-1$ we
conclude that
\bq
p_n^{(2)}(0)\ge [\PP(X_1<-\eps)]^k \overline{p}_n^{(2)}(y).\label{p0-py}
\eq
Of course, we also have the complementary trivial relations
$p_n^{(2)}(0)\le \overline{p}_n^{(2)} (0) \le p_n^{(2)}(y) \le
\overline{p}_n^{(2)}(y)$, so it suffices to consider
only $p_n^{(2)}(0)$ and $\overline{p}_n^{(2)}(0)$ which we
denote hereafter by $p_n^{(2)}$ and $\overline{p}_n^{(2)}$,
respectively. Obviously, $p_n^{(2)}$ and $\overline{p}_n^{(2)}$
have the same order (with $p_n^{(2)} = \overline{p}_n^{(2)}$
whenever $X_1$ has a density), and we consider both only
in order to draw the reader's attention to potential identities
connecting the two sequences $\{p_n^{(2)}\}$ and $\{\overline{p}_n^{(2)}\}$.

Persistence probabilities such as $p_n^{(2)}$ appear in many applications.
For example, the precise problem we consider here arises in the study
of the so-called sticky particle systems (c.f. \cite{Vysotsky2} and
the references therein). In case of standard normal $X_i$ it is also
related to entropic repulsion for $\nabla^2$-Gaussian fields (c.f.
\cite{grad2} and the references therein), though
here we consider the easiest version, namely a one dimensional
$\nabla^2$-Gaussian field. In his 1992 seminal paper, Sinai \cite{Sinai}
proved that if $\PP(X_1=1)=\PP(X_1=-1)=1/2$, then
$p_n^{(2)}\asymp n^{-1/4}$. However, his method relies on the fact
that for Bernoulli $\{X_k\}$ all local minima of $k \mapsto S_k^{(2)}$
correspond to values of $k$ where $S_k=0$, and as such form a
sequence of regeneration times. For this reason, Sinai's method
can not be readily extended to most other distributions. Using a
different approach, more recently Vysotsky \cite{Vysotsky} managed to
extend Sinai's result that $p_n^{(2)}\asymp n^{-1/4}$
to $X_i$ which are double-sided exponential, and a few other
special types of random walks. At about the same time,
Aurzada and Dereich \cite{AD} used strong approximation techniques
to prove the bounds $n^{-1/4}(\log n)^{-4}\lesssim p_n^{(2)}\lesssim
n^{-1/4}(\log n)^4$ for zero-mean random variables $\{X_i\}$
such that $\EE [e^{\beta |X_1|}] <\infty$ for some $\beta>0$.
However, even for $X_i$ which are standard normal variables it was not
known before whether these logarithmic terms are needed, and if not,
how to get rid of them. Our main result, stated below, resolves this
question,
at least when the lower tail of
$X_1^-=-\min(X_1,0)$ decays exponentially:\\
{\bf Decay Assumption:} We assume that there
exist constants $r>0$, $\theta> 1/r$ and finite $K$, $L$,
such that for all $t>0,s>0$,
\bq
\PP(-X_1> t+s)\le
K\PP(-X_1> t)\PP(-X_1>s) + L [\PP(-X_1>r)]^{\theta(t+s)}.
\label{decay}
\eq
\begin{remark}
With $\beta = - \limsup_{t \to \infty} \frac{1}{t} \log \PP(-X_1 > t)$, it
is easy to check that our condition holds (with $K=0$),
whenever $\PP(-X_1 > r) > e^{-\beta r}$ for some $r>0$.
Obviously this applies whenever $-X_1$ decays super-exponentially,
that is, for $\beta=\infty$. Conversely, considering $K=L>0$ and
$t=s=2^{\ell} r$, $\ell=0,1,\ldots$ it is not hard to show
that some exponential decay, namely $\beta>0$, is necessary
for our decay assumption. In the borderline case of an exponential
decay, that is, when $\beta>0$ is finite, our condition holds also
whenever $\PP(-X_1>t) \asymp  e^{-\beta t}$ (taking now $L=0$
and $K>0$, with $K=1$ corresponding to $X_1^-$ having a
\emph{New Better than Used} distribution).
\end{remark}
In this paper, we will prove:
\begin{theorem}\label{thm-main}
For i.i.d. $\{X_k\}$ of zero mean and $\EE|X_1|<\infty$,
let $S_n^{(2)}=S_1+S_2+\cdots +S_n$, where $S_i=X_1+X_2+\cdots+
X_i$. Then,
\bq
\sum_{k=0}^np_{k}^{(2)}\overline{p}_{n-k}^{(2)}\le c_1^2\frac{\EE|S_{n+1}|}{\EE|X_1|},\label{upper-bound}
\eq
where $c_1\le 6 \sqrt{30}$, and $c_1=2$ if $X_1$ is symmetric. Furthermore,
if $X_1^-$ satisfies the decay assumption (\ref{decay}), the
converse inequality
\bq
\sum_{k=0}^{n} p_{k}^{(2)}p_{n-k}^{(2)} \ge \frac{1}{c_2}
\frac{\EE|S_{n+1}|}{\EE|X_1|}
\label{lower-bound}
\eq
holds for some finite $c_2=c_2(K,L,\theta,r)$.
Taken together, these bounds imply that
\bq
\frac{1}{4c_1c_2}\sqrt{\frac{\EE|S_{n+1}|}{(n+1)\EE|X_1|}}\le p_n^{(2)}\le c_1\sqrt{\frac{\EE|S_{n+1}|}{(n+1)\EE|X_1|}}.\label{two-side}
\eq
\end{theorem}
\begin{remark}\label{Bernstein}
If $X_1$ has finite second moment, then by the \abbr{CLT},
there exist finite constants $C_2>C_1>0$ such that
$C_2 \sqrt{n} \ge \EE|S_n| \ge C_1 \sqrt{n}$ for all
$n \ge 1$. Consequently, we then have $p_n^{(2)}\asymp n^{-1/4}$ under
the decay assumption of Theorem~\ref{thm-main}. In contrast, for
any $0< \gamma < 1/4$ there exists integrable, zero-mean variable
$X_1$ for which $p_n^{(2)}\asymp n^{-\gamma}$. Indeed, considering
$\PP(Y_1 > y) = y^{-\alpha} 1_{y \ge 1}$ with $1 < \alpha <2$,
our decay assumption applies for the bounded below, zero-mean,
integrable random variable $X_1=Y_1-\EE Y_1
$. Setting $a_n=n^{1/\alpha}$,
clearly $n \PP(|X_1|>a_n x) \to x^{-\alpha}$ as $n \to \infty$,
hence $a_n^{-1} S_n - b_n$ converges in distribution to
a zero-mean, one-sided Stable$_\alpha$ variable $Z_\alpha$,
and it is further easy to check that
$b_n = a_n^{-1} n \EE[X_1 1_{|X_1| \le a_n}] \to b_\infty = -\EE Y_1$.
In fact, it is not hard to verify that
$\{a_n^{-1} S_n\}$ is a uniformly integrable sequence and
consequently $n^{-1/\alpha} \EE |S_n| \to \EE |Z_\alpha - \EE Y_1|$
finite and positive. From Theorem \ref{thm-main} we then
deduce that $p_n^{(2)}\asymp n^{-\gamma}$ for $\gamma=(1-1/\alpha)/2$.
\end{remark}
The sequences $\{S_k\}$ and $\{S_k^{(2)}\}$ are special cases
of the class of auto-regressive processes
$Y_k = \sum_{\ell=1}^L a_{\ell} Y_{k-\ell} + X_k$
with zero initial conditions, i.e. $Y_k \equiv 0$ when $k \le 0$
(where $S_k$ corresponds to $L=a_1=1$ and $S_k^{(2)}$ corresponds
to $L=a_1=2$, $a_2=-1$). While each of these stochastic processes
is a time-homogeneous Markov chain of state space $\reals$ and
$q_n = \PP(\tau>n)$ is merely the upper tail of the first
hitting time $\tau$ of $[y,\infty)$ by the relevant chain,
the general theory of Markov chains does not provide the
precise decay of $q_n$, which even in case $L=1$ ranges from
exponential decay for $a_1>0$ small enough, via the $O(n^{-1/2})$ decay
for $a=1$ to a constant $n \mapsto q_n$ in the
limit $a_1 \uparrow \infty$. While we do not pursue this here, we
believe that the approach we develop for proving Theorem \ref{thm-main}
can potentially determine the asymptotic behavior of $q_n$ for a
large collection of auto-regressive processes. This is of much
interest, since for example, as shown in \cite{Li-Shao},
the asymptotic tail probability that random Kac polynomials
have no (or few) real roots is determined in terms of
the limit as $r \to \infty$ of the power law tail decay
exponents for the iterates $S_k^{(r)} = \sum_{i=1}^k S_i^{(r-1)}$, $r \ge 3$.

Our approach further suggests that there might
be some identities connecting the sequences
$\{p_n^{(2)}\}$ and $\{\overline{p}_n^{(2)}\}$. Note that, if we denote
$$p_n^{(1)}=\PP\left(\max_{1\le k\le n}S_k<0\right),
\quad \overline{p}_n^{(1)}=\PP\left(\max_{1\le k\le n}S_k\le 0\right),$$
then as we show in Section \ref{pn1-proof}
there are indeed identities connecting the sequences
$\{p_n^{(1)}\}$ and $\{\overline{p}_n^{(1)}\}$.
As we mentioned earlier, it is well-known that $p_n^{(1)}$ is
of the order $n^{-1/2}$. The next proposition provides
the exact value of $p_n^{(1)}$ for symmetric random variables
with a density,
and the elegant argument of its proof
serves as the starting point of our approach
to the study of $p_n^{(2)}$.
\begin{proposition}\label{pn1}
If $X_i$ are mean zero i.i.d. symmetric random variables
then for all $n \geq 1$,
 \bq
p_n^{(1)} \le \frac{(2n-1)!!}{(2n)!!}\le \overline{p}_n^{(1)}\,.
\label{discrete}
\eq
In particular, if $X_1$ also has a density, then
\bq
p_n^{(1)} = \frac{(2n-1)!!}{(2n)!!}\,.\label{p1}\eq
\end{proposition}

\begin{remark}
Let $B(s)$ denote a Brownian motion starting at $B(0)=0$ and
consider the integrated Brownian motion $Y(t)=\int_0^t B(s) ds$.
Sinai \cite{Sinai} proved the existence of positive constants
$A_1$ and $A_2$ such that for any $T>0$,
\begin{equation}\label{eq:sinai-bm}
A_1 T^{-1/4}\le \PP\Big(\sup_{t\in [0,T]} Y(t) \le 1\Big)
\le A_2 T^{-1/4}.
\end{equation}
Upon setting $\eps = T^{-3/2}$ and $t=uT$,
by Brownian motion scaling this is equivalent up to a
constant to the following result that can be derived from an
implicit formula of McKean \cite{McKean} (c.f. \cite{Simon}):
$$
\lim_{\eps\to 0^+}\eps^{-1/6}\PP\Big(\sup_{u\in [0,1]} Y(u)
\le \eps\Big)=\frac{3\Gamma(5/4)}{4\pi\sqrt{2\sqrt{2\pi}}}.
$$
Since the iterated partial sums $S_n^{(2)}$
corresponding to i.i.d standard normal random
variables $\{X_i\}$, forms a discretization of $Y(t)$,
the right-most inequality in (\ref{eq:sinai-bm})
readily follows from Theorem \ref{thm-main}.
Indeed, with $\EE [Y(k)Y(m)]=k^2 (3m-k)/6$ and
$\EE[S_k^{(2)} S_m^{(2)}]=k(k+1)(3m-k+1)/6$ for
$m \ge k$, setting $Z(k)=\sqrt{(1+1/k)(1+1/(2k)}Y(k)$,
results with $\EE [(S_k^{(2)})^2]=\EE[Z(k)^2]$ and it
is further not hard to show that
$f(m,k):= \EE [S_m^{(2)}S_k^{(2)}]/\EE [Z(m)Z(k)] \ge 1$
for all $m \ne k$ (as $f(k+1,k) \ge 1$ and $d f(m,k)/dm >0$
for any $m \ge k+1$). Thus, by Slepian's lemma, we have that for any $y$
$$
\PP\Big(\max_{1\le k\le n} Z(k) <  y\Big)\le p_n^{(2)}(y) \,,
$$
and setting $n$ as the integer part of $T\ge 1$ it follows that
$$
\PP\Big(\sup_{t\in [0,T]} Y(t) \le 1 \Big)\le
\PP\Big(\max_{1\le k\le n}Y(k)\le 1\Big)
\le \PP\Big(\max_{1\le k\le n}Z(k) < 2 \Big) \le p_n^{(2)}(2) \,.
$$
Since $p_n^{(2)}(2) \le c p_n^{(2)}$ for some finite constant $c$ and all
$n$, we conclude from Theorem \ref{thm-main} that
$$
\PP\Big(\sup_{t\in [0,T]} Y(t) \le 1\Big) \le 2 c (n+1)^{-1/4}\le 2 c T^{-1/4}\,.
$$
\end{remark}

\section{Proof of Proposition \ref{pn1}}\label{pn1-proof}
Setting $S_0 =0$ let $M_n=\max_{0 \le j \le n} S_j$ and
consider the $\{0, 1,2,\ldots,n\}$-valued random variable
$$
\cN = \min\left\{ l \ge 0 : S_l=M_n \right\} \,.
$$
For each $k=1,2,\ldots,n-1$ we have that
\begin{align}
\{\cN=k\}&=\{X_k>0, X_{k}+X_{k-1}>0,\ldots,
X_k+X_{k-1}+\cdots+X_1>0;
\nonumber \\
& X_{k+1}\le 0, X_{k+1}+X_{k+2}\le
0, \ldots, X_{k+1}+X_{k+2}+\cdots+X_n\le 0\}.
\nonumber
\end{align}
By the independence of $\{X_i\}$, the latter identity implies that
\begin{align}
\P(\cN=k)&=\P(X_k>0, X_{k}+X_{k-1}>0,\ldots, X_k+X_{k-1}+\cdots+X_1>0)
\nonumber \\
\nonumber
&\times
\P(X_{k+1}\le 0, X_{k+1}+X_{k+2}\le 0, \ldots, X_{k+1}+X_{k+2}+\cdots+X_n\le 0)\\
&=p_k^{(1)}\overline{p}_{n-k}^{(1)},
\nonumber
\end{align}
where the last equality follows
from our assumptions that $X_i$ are i.i.d symmetric random
variables.  Also note that $\PP(\cN=0)=\overline{p}_n^{(1)}$
and
$$\PP(\cN=n)=\P(X_n>0,
X_{n}+X_{n-1}>0,\ldots, X_n+X_{n-1}+\cdots+X_1>0)=p_n^{(1)}.$$
Thus, setting $p_0^{(1)}=\overline{p}_0^{(1)}=1$ we arrive at
the identity
\bq
\sum_{k=0}^{n}p_{k}^{(1)}\overline{p}_{n-k}^{(1)}=\sum_{k=0}^n\P(\cN=k)=1,\label{identity}
\eq
holding for all $n\ge 0$.

Fixing $x \in [0,1)$, upon multiplying (\ref{identity})
by $x^n$ and summing over $n \ge 0$, we arrive
at $P(x) \overline{P}(x) = \frac1{1-x}$, where
$P(x)=\sum_{k=0}^\infty p_k^{(1)}x^k$ and
$\overline{P}(x)=\sum_{k=0}^\infty \overline{p}_k^{(1)}x^k$.
Now, if $X_1$ also has a density then $p_k^{(1)}=\overline{p}_k^{(1)}$ for all $k$
and so by the preceding $P(x)=\overline{P}(x)=(1-x)^{-1/2}$. Consequently,
$p_n^{(1)}$ is merely the coefficient of $x^n$ in the Taylor
expansion at $x=0$ of the function $(1-x)^{-1/2}$, from which
we immediately deduce the identity (\ref{p1}).

If $X_1$ does not have a density, let $\{Y_i\}$ be i.i.d. standard normal random variables, independent of the sequence
$\{X_i\}$ and denote by $S_k$ and $\tilde{S}_k$ the partial sums of $\{X_i\}$ and $\{Y_i\}$, respectively.
Note that for any $\eps>0$, each of the i.i.d. variables $X_i+\eps Y_i$ is symmetric and has a density,
with the corresponding partial sums being $S_k+\eps \tilde{S}_k$. Hence, for any $\delta>0$ we have that
\begin{align}
\PP\Big(\max_{1\le k\le n} S_k<-\delta\Big)&\le
\PP\Big(\max_{1\le k\le n} (S_k+\eps \tilde{S}_k)\le 0\Big)+\PP\Big(\max_{1\le k\le n} \eps\tilde{S}_k \ge\delta\Big)
\nonumber\\
&=\frac{(2n-1)!!}{(2n)!!}+\PP\Big(\max_{1\le k\le n} \eps\tilde{S}_k \ge\delta\Big).
\nonumber
\end{align}
Taking first $\eps\downarrow 0$ followed by $\delta\downarrow 0$, we conclude that
$$
\PP\Big(\max_{1\le k\le n} S_k<0\Big)\le \frac{(2n-1)!!}{(2n)!!}\,,
$$
and a similar argument works for the remaining inequality in (\ref{discrete}).

\begin{remark} The argument of the next section allows us to modify this proof and deduce order of
$(n+1)^{-1/2}$ upper and lower bounds for $p_n^{(1)}$ even in the non-symmetric case. However, since
this result is already known, we shall not do so here.
\end{remark}

\section{Proof of Theorem \ref{thm-main}}\label{main-proof}

By otherwise considering $X_i/\EE|X_i|$, we assume without loss of
generality that $\EE|X_1|=1$.
To adapt the method of Section \ref{pn1-proof} for dealing with
the iterated partial sums $S_n^{(2)}$, we introduce the parameter $t\in \RR$ and
consider the iterates $S_j^{(2)}(t) =S_0(t) +\cdots+S_j(t)$, $j \ge 0$,
of the translated partial sums $S_k(t) =t+S_k^{(1)}$, $k \ge 0$.
That is, $S_j^{(2)}(t)=(j+1) t + S_j^{(2)}$ for each $j \ge 0$.

Having fixed the value of $t$, we define the following
$\{0,1,2,\ldots,n\}$-valued random variable
$$
\cK_t=\min\Big\{l \ge 0: S_l^{(2)}(t) =\max_{0\le j\le n} S_j^{(2)}(t)\Big\}.
$$
Then, for each $k=2,3,\ldots,n-2$, we have the identity
\begin{align}
\{\cK_t=k\}&=\left\{S_k(t) >0,
S_{k}(t) +S_{k-1}(t) >0,\ldots,
S_k(t) +S_{k-1}(t) +\cdots+S_1(t) >0;\right.
\nonumber
\\&\ \
\left. S_{k+1}(t) \le 0, S_{k+1}(t) +S_{k+2}(t) \le 0, \ldots,
S_{k+1}(t) +S_{k+2}(t) +\cdots+S_n(t) \le 0\right\}
\nonumber \\
&= \left\{S_{k}(t)>0; X_{k}<2S_k(t),\ldots,(k-1)X_k+\cdots+X_2<kS_k(t)\right\} \cap
\{S_{k+1}(t)\le 0\}
\nonumber
\\&\cap\left\{X_{k+2}\le -2S_{k+1}(t), \ldots,
(n-k-1)X_{k+2}+\cdots +X_n\le -(n-k)S_{k+1}(t)\right\}.
\nonumber
\end{align}
Next, for $2 \le k \le n$ we define $Y_{k,2} \in \sigma(X_2,\ldots,X_k)$ and
$Y_{k,n} \in \sigma(X_k,\ldots,X_n)$ such that
\begin{align}
\nonumber
&Y_{k,2}=\max\left\{\frac{X_k}{2},\frac{2X_k+X_{k-1}}{3}, \ldots, \frac{(k-1)X_{k}+\cdots+X_2}{k}\right\},\\
&
Y_{k,n}=\max\left\{\frac{X_k}{2},\frac{2X_k+X_{k+1}}{3}, \ldots, \frac{(n-k+1)X_k+\cdots+X_n}{n-k+2}\right\}.
\nonumber
\end{align}
It is then not hard to verify that the preceding identities translate into
\begin{align}
\{\cK_t=k\}&= \{S_k(t)>0 \ge S_{k+1}(t) \} \cap \{Y_{k,2}<S_k(t)\}\cap\{Y_{k+2,n}\le -S_{k+1}(t)\} \label{ID} \\
&= \{ -S_k + (Y_{k,2})^+ < t \le -X_{k+1} - S_k - (Y_{k+2,n})^+ \}
\label{ID2}
\end{align}
holding for each $k=2,\ldots,n-2$. Further, for $k=1$ and $k=n-1$ we have that
\begin{align}
\{\cK_t=1\}  &=\{S_1(t)>0\}\cap \{S_{2}(t)\le 0\}\cap\{Y_{3,n}\le -S_{2}(t)\}\,,\nonumber \\
\{\cK_t=n-1\}&=\{S_{n-1}(t)>0\}\cap\{Y_{n-1,2}<S_{n-1}(t)\}\cap \{S_{n}(t)\le 0\}\,,\nonumber
\end{align}
so upon setting $Y_{1,2}=Y_{n+1,n}=-\infty$, the identities (\ref{ID}) and (\ref{ID2})
extend to all $1 \le k \le n-1$.

For the remaining cases, that is, for $k=0$ and $k=n$, we have instead that
\begin{align}
\{\cK_t=0\}&=\{t\le -X_1 - (Y_{2,n})^+ \}\,,\label{0}\\
\{\cK_t=n\}&=\{-S_n+(Y_{n,2})^+ < t\}\,.\label{n}
\end{align}

In contrast with the proof of Proposition \ref{pn1}, here we have
events $\{(Y_{k,2})^+<S_k(t)\}$ and $\{(Y_{k+2,n})^+\le -S_{k+1}(t)\}$
that are linked through $S_{k}(t)$ and consequently
not independent of each other. Our goal is to unhook
this relation and in fact the parameter $t$
was introduced precisely for this purpose.

\subsection{Upper bound}
For any integer $n>1$, let
$$A_n =\max_{1\le k\le n}\{-S_{k+1} \},\ \ B_n =-\max_{1\le k\le n}\{S_k \}.$$
By definition $A_n\ge B_n$. Further, for any $1 \le k \le n-1$, from (\ref{ID}) we have that the event $\{\cK_t=k\}$ implies that
$\{S_k(t) > 0 \ge S_{k+1}(t) \} = \{ -S_k < t \leq - S_{k+1} \}$
and hence that $\{B_{n-1} <t \le A_{n-1} \}$. From (\ref{ID2}) we also
see that for any $1 \le k \le n-1$,
$$
\int_{\RR} 1_{\{\cK_t=k\}}dt \ge
(X_{k+1})^- 1_{\{Y_{k,2}<0\}}1_{\{Y_{k+2,n}\le 0\}}
$$
and consequently,
\begin{align}
A_{n-1}-B_{n-1} = \int_{\RR} 1_{\{B_{n-1} <t \le A_{n-1} \}}dt &\ge
\sum_{k=1}^{n-1} \int_{\RR} 1_{\{\cK_t=k\}}dt \label{sum-indicator}\\
&\ge \sum_{k=1}^{n-1} (X_{k+1})^- \, 1_{\{Y_{k,2}<0\}}1_{\{Y_{k+2,n}\le 0\}}\,. \nonumber
\end{align}
Taking the expectation of both sides
we deduce from the mutual independence of $Y_{k,2}$, $X_{k+1}$ and
$Y_{k+2,n}$ that
$$
\EE [A_{n-1}-B_{n-1}] \ge \sum_{k=1}^{n-1}\EE [(X_{k+1})^-] \PP(Y_{k,2}<0) \PP(Y_{k+2,n}\le 0) \,.
$$
Next, observe that since the sequence $\{X_i\}$ has an exchangeable law,
\begin{align}
\PP(Y_{k,2} < 0) &= \PP(X_k < 0, 2 X_k + X_{k-1} < 0, \ldots, (k-1) X_k + \cdots + X_2 < 0)
\nonumber \\
&= \PP(X_1<0, 2 X_1 + X_2 < 0, \ldots, (k-1) X_1 + \cdots + X_{k-1} < 0) = p_{k-1}^{(2)} \,.
\label{basic-bd}
\end{align}
Similarly, $\PP(Y_{k+2,n} \le 0) = \overline{p}_{n-1-k}^{(2)}$.
With $X_{k+1}$ having
zero mean, we have that $\EE [(X_{k+1})^-] = \EE[ (X_{k+1})^+] = 1/2$ (by
our assumption that $\EE |X_{k+1}| = \EE |X_1| =1$). Consequently,
for any $n>2$,
$$
\EE [A_{n-1}-B_{n-1}] \ge \frac{1}{2} \sum_{k=1}^{n-1} p_{k-1}^{(2)}\overline{p}_{n-1-k}^{(2)}
= \frac{1}{2} \sum_{k=0}^{n-2} p_k^{(2)} \overline{p}_{n-2-k}^{(2)} \,.
$$
With $\EE[S_{n+1}]=0$ and $\{X_k\}$ exchangeable, we clearly have that
\begin{equation}\label{an-bn:bd}
\EE[A_n-B_n] = \EE[\max_{1 \le k \le n} S_k]
+ \EE[\max_{1 \le k \le n} \{S_{n+1}-S_{k+1}\}]
=2 \EE[\max_{1 \le k \le n} S_k] \,.
\end{equation}
Recall Ottaviani's maximal inequality that for a symmetric random walk
$\PP(\max_{k=1}^n S_k \ge t) \le 2 \PP(S_n \ge t)$ for any $n, t \ge 0$,
hence in this case
$$
\EE[\max_{1 \le k \le n} S_k] \le
2 \int_0^\infty \PP(S_n\ge t) dt = \EE |S_n| \,.
$$
To deal with the general case, we replace Ottaviani's maximal inequality
by Montgomery-Smith's inequality
$$
\PP(\max_{1\le k\le n}|S_k|\ge t)
\le 3\max_{1\le k\le n}\PP(|S_k|\ge t/3)\le 9\PP(|S_n| \ge t/30)
$$
(see \cite{Montgomery}), from which we deduce that
\begin{equation}\label{an-bn:bd2}
\EE[\max_{1 \le k \le n} S_k] \le
9 \int_0^\infty \PP(|S_n|\ge t/30) dt = 270 \EE |S_n|
\end{equation}
and thereby get (\ref{upper-bound}).
Finally, since $n \mapsto p_n^{(2)}$ is non-increasing and
$p_n^{(2)}\le\overline{p}_{n}^{(2)}$, the upper bound of
(\ref{two-side}) is an immediate consequence of (\ref{upper-bound}).

\subsection{Lower bound}
Turning to obtain the lower bound, let
$$m_n:= -X_1-(Y_{2,n})^+\,,\ \ M_n:=-S_n + (Y_{n,2})^+ \,.$$
Note that for any $n \ge 2$,
$$
Y_{n,2}+Y_{2,n} \ge \frac{1}{n} [(n-1) X_n + \cdots + X_2] + \frac{1}{n} [
(n-1) X_2 + \cdots + X_n] = S_n - X_1 \,,
$$
and consequently,
\begin{equation}\label{lbd-diff}
M_n-m_n \ge X_1-S_n + (Y_{2,n}+Y_{n,2})^+ \ge (X_1-S_n)^+ = (X_2+\cdots+X_n)^-\,.
\end{equation}
In particular, $M_n\ge m_n$. From (\ref{0}) and (\ref{n}) we know that if $m_n<t \le M_n$
then necessarily $1 \le \cK_t \le n-1$. Therefore,
\begin{equation}\label{L}
M_n-m_n = \int_{\RR} 1_{\{m_n < t \le M_n\}}dt \le
\sum_{k=1}^{n-1} \int_{\RR} 1_{\{\cK_t=k\}} dt \,.
\end{equation}
In view of (\ref{ID2}) we have that for any $1 \le k \le n-1$,
$$
b_k := \EE [\int_{\RR}1_{\{\cK_t=k\}}dt]
=\EE \Big[ \left(X_{k+1} + (Y_{k,2})^+ \, +(Y_{k+2,n})^+ \right)^-\Big].
$$
By the mutual independence of the three variables on the right side,
denoting by $F_{k,2}$ and $F_{k+2,n}$ the distribution functions of $(Y_{k,2})^+$ and
$(Y_{k+2,n})^+$, respectively, we thus find that
\begin{align}
b_k &=\int_0^\infty\int_0^\infty \EE[(X_{k+1}+x+y)^-] dF_{k,2}(x)\, dF_{k+2,n}(y)\nonumber\\
&=\int_0^\infty\int_0^\infty g_{k+1} (x+y) \,dF_{k,2}(x)\, dF_{k+2,n}(y)\,, \nonumber
\end{align}
where $g_k (t)=\int_0^\infty\PP(-X_k>t+u)du$.
Since $\{X_k\}$ have identical distributions, $g_k(t)=g_1 (t)$ does not
depend on $k$ and we have already seen that $g_1(0)=\EE[(X_1)^-]=1/2$.
Thus, setting $\alpha=-\log\PP(-X_1>r)$,
applying our decay assumption (\ref{decay}) first for
$t=x+y$, $s=u$ and then for $t=x$, $s=y$, we find that
\begin{align}
g_1 (x+y)&\le K \PP(-X_1>x+y) g_1(0) +
L \int_0^\infty e^{-\theta \alpha (x+y+u)} du \nonumber \\
&\le \frac{K^2}{2} \PP(-X_1>x)\PP(-X_1>y) + L_1
e^{-\theta\alpha(x+y)} \, \nonumber
\end{align}
where $L_1 = L [K/2 + (\theta \alpha)^{-1} ]$.
Consequently, we arrive at the bound
\begin{align}
b_k &\le \frac{K^2}{2}
 \int_0^\infty \PP(-X_1>x)
 dF_{k,2}(x) \int_0^\infty \PP(-X_1 > y) dF_{k+2,n}(y)
\nonumber \\
&+ L_1 \int_0^\infty e^{-\theta\alpha x} dF_{k,2} (x)
dF_{k+2,n}(y)\nonumber \\
&=
\frac{K^2}{2}\PP(X_{k+1} + (Y_{k,2})^+ <0) \PP(X_{k+1}+(Y_{k+2,n})^+ < 0)
+ L_1 \EE[ e^{-\theta \alpha (Y_{k,2})^+} ]
\EE[ e^{-\theta \alpha (Y_{k+2,n})^+} ] \,. \nonumber
\end{align}
Next, observe that
just as we did in deriving the identity (\ref{basic-bd}),
\begin{align}
\PP(X_{k+1}+(Y_{k,2})^+  < 0) &= \PP(
X_{k+1}<0, 2 X_{k+1} + X_k < 0, \ldots, k X_{k+1} + \cdots + X_2 < 0)
\nonumber \\
&= \PP(X_1<0, 2 X_1 + X_2 < 0, \ldots, k X_1 + \cdots + X_k < 0) = p_k^{(2)} \,.
\nonumber
\end{align}
By a similar reasoning one verifies that
$p_{n-k}^{(2)}= \PP(X_{k+1}+(Y_{k+2,n})^+<0)$.
Now, by exchangeability of $\{X_k\}$,
the sequence $k \mapsto \PP(Y_{k,2} <z)$ is non-increasing
for any fixed value of $z$. Furthermore, if $z \ge r > 0$ then
\begin{align}
\PP(Y_{k,2}<z-r) & =
\PP(X_k < 2 (z-r),
\ldots, (k-1)X_{k}+\cdots+X_2< k (z-r))
\nonumber \\
&\ge \PP(X_k<-r)\PP\big(
X_{k-1}<2z, \ldots, (k-2)X_{k-1}+\cdots+X_2<(k-1) z \big) \nonumber \\
&=\PP(-X_1 > r) \PP(Y_{k-1,2}<z) \ge \PP(X_1<-r) \PP(Y_{k,2}<z). \nonumber
\end{align}
Iterating this inequality for $z=j r$, $j=1,\ldots,k:=\ceil{x/r}$
we deduce that for any $x>0$,
$$
\PP(Y_{k,2}<x) \le \PP(Y_{k,2}<k r) \le [\PP(-X_1>r)]^{-k} \PP(Y_{k,2}<0)
\le e^{\alpha (x/r+1)} p_{k-1}^{(2)}
$$
(relying on the identity (\ref{basic-bd}) for the right-most inequality).
Consequently,
$$
\EE[ e^{-\theta \alpha (Y_{k,2})^+} ] =
\theta \alpha \int_0^\infty \PP(0 \le Y_{k,2} < x) e^{-\theta \alpha x} dx
\le \kappa p_{k-1}^{(2)} \,,
$$
with $\kappa = e^{\alpha} \theta/(\theta-1/r)$ finite.
Similarly, we find that for any $y>0$,
$$
\PP(Y_{k+2,n}<y)=p_{n-k-1}^{(2)}(y) \le e^{\alpha (y/r+1)} p_{n-k-1}^{(2)},
$$
hence
$$
\EE[ e^{-\theta \alpha (Y_{k+2,n})^+} ] \le \kappa p_{n-k-1}^{(2)} \,.
$$
Combining all these bounds we have by the monotonicite of $k \mapsto
p^{(2)}_k$ that
$$
b_k \le \frac{K^2}{2} p_k^{(2)} p_{n-k}^{(2)} +
L_1 \kappa^2 p_{k-1}^{(2)} p_{n-k-1}^{(2)}
\le \frac{c_2}{2} p_{k-1}^{(2)} p_{n-k-1}^{(2)} \,,
$$
where $c_2 = K^2 + 2 L_1 \kappa^2$ is a finite constant.
Thus, considering the expectation of both sides of (\ref{L})
we deduce that for any $n>2$,
\begin{eqnarray*}
\EE(M_n-m_n) \le \frac{c_2}{2}
\sum_{k=1}^{n-1}p_{k-1}^{(2)}p_{n-k-1}^{(2)}\,.
\end{eqnarray*}
In view of (\ref{lbd-diff}) we also have that
$\EE(M_n-m_n)\ge \EE [(S_{n-1})^-] = \frac12 \EE|S_{n-1}|$, from which
we conclude that (\ref{lower-bound}) holds for all $n \ge 1$.

Turning to lower bound $p_n^{(2)}$ as in (\ref{two-side}),
recall that $n \mapsto p_n^{(2)}$ is non-increasing. Hence,
applying (\ref{lower-bound}) for $n=2m+1$ and utilizing the
previously derived upper bound of (\ref{two-side}) we have that
\begin{align}
\frac1{c_2} \EE|S_{2(m+1)}|&\le
2\sum_{k=0}^mp_k^{(2)}p_m^{(2)}\le 2c_1p_m^{(2)}
\sum_{k=0}^m\sqrt{\frac{\EE|S_{k+1}|}{k+1}}\nonumber
\\& \le 4 c_1 p_m^{(2)}\sqrt{(m+1)\EE|S_{m+1}|}\,,
\end{align}
where in the last inequality we use the fact that
for independent, zero-mean $\{X_k\}$, the sequence
$|S_k|$ is a sub-martingle, hence
$k \mapsto \EE|S_k|$ is non-decreasing.
This proves the lower bound of (\ref{two-side}).

\end{document}